\documentclass[11pt,a4paper]{article}

\usepackage[utf8]{inputenc}
\usepackage[T1]{fontenc}
\usepackage{amsmath,amssymb,amsthm}
\usepackage{algorithm}
\usepackage{algpseudocode}
\usepackage{graphicx}
\usepackage{xcolor}
\usepackage{hyperref}
\usepackage{geometry}
\usepackage{natbib}
\usepackage{listings}
\usepackage{booktabs}

\hypersetup{
    colorlinks=true,
    linkcolor=blue,
    citecolor=blue,
    urlcolor=blue
}

\newtheorem{definition}{Definition}

\newtheorem{proposition}{Proposition}

\title{Evaluating Numerical Accuracy in Mixed-Precision Computing by Dual-Delta Testing}
\author{Peichen Xie}
\date{}

\begin{document}

\maketitle

\begin{abstract}
Mixed-precision computing has become increasingly important in modern high-performance computing and machine learning applications. When implementing custom mixed-precision functions---such as fused operators, optimized GPU kernels, or quantized inference paths---it is critical to verify their numerical accuracy. Traditional approaches typically compare the custom implementation against a reference using a single error metric. However, this single-delta approach provides limited insight into whether the observed errors are inherent to the precision level or specific to the implementation. This paper introduces \textit{Dual-Delta Testing}, a systematic methodology that evaluates two error distributions against a high-precision oracle, enabling rigorous comparison between a custom implementation and a baseline reference. We present the mathematical framework, algorithmic formulation, statistical analysis techniques, and practical examples demonstrating the methodology's effectiveness in evaluating numerical accuracy.
\end{abstract}

\section{Introduction}

\subsection{Motivation}

The proliferation of mixed-precision arithmetic in modern computing has fundamentally transformed how we approach numerical computation across diverse application domains. This transformation is driven by a compelling convergence of performance, efficiency, and hardware innovation. Lower-precision arithmetic operations, such as 16-bit floating-point (FP16) and brain-float (BF16) formats, execute significantly faster than their full-precision counterparts while maintaining precision sufficient for algorithm~\citep{micikevicius_mixed_2018,kalamkar_study_2019}. This performance advantage becomes particularly pronounced in data-intensive applications where billions of operations must be performed in real-time.

Beyond raw computational speed, reduced precision data types offer substantial memory efficiency gains. By halving or quartering the memory footprint compared to traditional 32-bit or 64-bit representations, these formats enable practitioners to train larger neural networks, process bigger datasets, and deploy models on resource-constrained devices~\citep{gupta_deep_2015,wang_training_2018}. Modern hardware accelerators, such as NVIDIA's Tensor Cores and Google's Tensor Processing Units (TPUs), have been specifically designed for lower-precision operations, providing order-of-magnitude speedups for mixed-precision workloads~\citep{jouppi_-datacenter_2017}.

However, this shift toward reduced precision introduces a critical challenge that every practitioner must confront: numerical accuracy. As we compress our numeric representations, we inevitably introduce approximation errors that can propagate through complex computational pipelines. When developing custom mixed-precision implementations---whether optimizing GPU kernels for matrix operations, designing fused operators for neural network inference, or implementing quantized arithmetic for edge deployment---engineers and researchers face a fundamental and often vexing question: \textit{Is my implementation numerically correct, and how does it compare to existing baseline methods?}

\subsection{The Problem with Single-Delta Testing}

The conventional approach to validating numerical accuracy involves computing a single error metric (delta) between the custom implementation and a reference implementation:
\begin{equation}
    \Delta_{\text{single}} = \{\text{Error}(\text{impl}_{\text{custom}}, \text{impl}_{\text{ref}})\}_{\text{inputs}}
\end{equation}

While this single-delta methodology appears straightforward, it suffers from conceptual and practical limitations that can lead to incorrect conclusions about implementation correctness. The primary issue is the lack of contextual information: when we observe a large error between two implementations, we cannot definitively determine whether this discrepancy arises from fundamental precision limitations shared by both implementations or from a specific bug in the custom implementation under test. This ambiguity becomes particularly problematic in mixed-precision computing, where numerical errors are expected but their magnitude and distribution must be carefully characterized.

Further complicating matters is the question of \emph{comparative} accuracy: even if a custom mixed-precision implementation appears ``close'' to a reference, what does that imply? A small discrepancy does not necessarily mean the custom implementation is correct, nor does a larger discrepancy necessarily indicate a bug---both results may simply reflect different rounding, accumulation order, fused operations, or hardware-specific behaviors. Moreover, the comparison itself is ambiguous: should we conclude that the custom implementation is \emph{more accurate} than the reference, \emph{equally accurate}, or \emph{less accurate}? Answering that requires more than a single error number. We need a principled way to decide whether one implementation exhibits systematically larger errors.

\subsection{The Solution of Dual-Delta Testing}

Dual-Delta Testing addresses these limitations through a paradigm shift: rather than directly comparing two implementations of mixed-precision arithmetic, we introduce a high-precision oracle that serves as a trusted ground truth, and we measure two separate error distributions relative to this oracle:
\begin{align}
    \Delta_1 &= \{\text{Error}(\text{impl}_1, \text{oracle})\}_{\text{inputs}} \\
    \Delta_2 &= \{\text{Error}(\text{impl}_2, \text{oracle})\}_{\text{inputs}}
\end{align}

This dual-distribution framework transforms the validation problem from a single metric into a rich statistical analysis. By evaluating both implementations against the same high-precision baseline, we can definitively determine whether the custom implementation achieves comparable accuracy to the reference implementation. The distributional nature of the approach reveals whether errors are systematic (indicating potential algorithmic issues) or random (suggesting inherent precision limitations). We can precisely quantify the magnitude of precision-related degradation relative to the mathematical ground truth, and make statistically rigorous statements about numerical correctness that would be impossible with a single error metric.

\section{Mathematical Framework}

\subsection{Notation and Definitions}

\begin{definition}[Implementation]
An implementation is a function $f: \mathcal{X} \rightarrow \mathcal{Y}$ that maps inputs from domain $\mathcal{X}$ to outputs in codomain $\mathcal{Y}$, computed using a specific numerical precision and algorithm.
\end{definition}

\begin{definition}[Oracle]
An oracle $f_\Omega: \mathcal{X} \rightarrow \mathcal{Y}$ is a high-precision implementation that serves as the ground truth. Typically, the oracle uses significantly higher precision (e.g., FP64) than the implementations being tested.
\end{definition}

\begin{definition}[Error Metric]
An error metric is a function $\epsilon: \mathcal{Y} \times \mathcal{Y} \rightarrow \mathbb{R}^+$ that quantifies the discrepancy between two outputs.
\end{definition}

\subsection{The Dual-Delta Formulation}

Given two implementations ($f_1$ as custom and $f_2$ as reference), an oracle $f_\Omega$, an error metric $\epsilon$, and an input distribution $P(x)$ over domain $\mathcal{X}$, we formalize the dual-delta methodology. For a sample of $N$ inputs $\{x_1, x_2, \ldots, x_N\}$ drawn from $P(x)$, we compute:

\begin{align}
    \Delta_1 &= \{\epsilon(f_1(x_i), f_\Omega(x_i))\}_{i=1}^N \\
    \Delta_2 &= \{\epsilon(f_2(x_i), f_\Omega(x_i))\}_{i=1}^N
\end{align}

The dual-delta testing evaluates the relationship between these distributions through statistical comparison. 

\begin{proposition}[Equivalence]
If the distributions $\Delta_1$ and $\Delta_2$ are statistically indistinguishable, we conclude that $f_1$ and $f_2$ exhibit comparable numerical accuracy relative to the oracle.
\end{proposition}

\begin{proposition}[Numerical Accuracy]
If $\Delta_1$ is statistically significantly smaller than $\Delta_2$, we conclude that $f_1$ is numerically more accurate than $f_2$. Conversely, if $\Delta_1$ is significantly larger, we conclude that $f_2$ is more accurate.
\end{proposition}

\begin{proposition}[Numerical Stability]
If the variance of $\Delta_1$ is statistically significantly smaller than the variance of $\Delta_2$, we conclude that $f_1$ is numerically more stable than $f_2$. Conversely, the variance of $\Delta_1$ is significantly larger, we conclude that $f_2$ is more stable.
\end{proposition}

\section{Algorithmic Formulation}

\subsection{The Core Algorithm}

Algorithm~\ref{alg:dual_delta} presents the pseudocode for a Dual-Delta Test.

\begin{algorithm}
\caption{Dual-Delta Test}
\label{alg:dual_delta}
\begin{algorithmic}[1]
\Require $f_1, f_2$: implementations to compare
\Require $f_\Omega$: oracle implementation
\Require $\textsc{GenerateInput}()$: input generator
\Require $\epsilon$: error metric
\Require $N$: number of test cases
\Ensure $\Delta_1, \Delta_2$: error distributions
\For{$i = 1$ to $N$}
    \State $x \gets \textsc{GenerateInput}()$
    \State $y_1 \gets f_1(x)$
    \State $y_2 \gets f_2(x)$
    \State $y_\Omega \gets f_\Omega(x)$
    \State $\Delta_1.\textsc{Append}(\epsilon(y_1, y_\Omega))$
    \State $\Delta_2.\textsc{Append}(\epsilon(y_2, y_\Omega))$
\EndFor
\State \Return $\Delta_1, \Delta_2$
\end{algorithmic}
\end{algorithm}

\subsection{Implementation Considerations}

\subsubsection{Input Generation}
The design of the input generator $\textsc{GenerateInput}()$ critically influences the validity and generalizability of test results. Random sampling from appropriate statistical distributions (such as Gaussian or uniform distributions) provides broad coverage of the input space while maintaining computational efficiency. However, random sampling alone proves insufficient---edge cases that occur rarely in random sampling often trigger the most severe numerical errors. Domain-specific patterns add another critical dimension: real-world data rarely follows idealized statistical distributions, so incorporating realistic data characteristics (such as sparsity patterns in scientific computing or typical activation distributions in neural networks) ensures that test results meaningfully predict production behavior.

\subsubsection{Oracle Selection}
The oracle must satisfy a fundamental precision separation property: its errors relative to the mathematically exact result must be negligibly small compared to the errors of the implementations under test. Formally:
\begin{equation}
    \epsilon(f_\Omega(x), f_{\text{true}}(x)) \ll \epsilon(f_i(x), f_{\text{true}}(x)) \quad \forall i \in \{1,2\}
\end{equation}
where $f_{\text{true}}$ represents the mathematically exact result. In practice, this requirement translates to using significantly higher precision for the oracle than the implementations being tested. For single-precision (FP32) implementations, double-precision (FP64) arithmetic typically provides adequate separation. For lower-precision (such as FP16, BF16, and FP8) implementations, both FP64 and FP32 arithmetic can suffice. For applications demanding even higher confidence---such as financial calculations or safety-critical systems---arbitrary precision libraries like MPFR \cite{fousse_mpfr_2007} eliminate oracle error entirely at the cost of computational overhead.

\subsubsection{Error Metric Selection}
The choice of error metric profoundly shapes our interpretation of numerical accuracy, and no single metric universally dominates across all applications. Among the most commonly used error metrics, we recommend considering the following options:

\begin{enumerate}
    \item Norm-wise relative error $\epsilon_{\text{rel}}(y, y') = \frac{\|y - y'\|}{\|y'\|}$ measures the overall relative difference between outputs, making it invaluable when outputs span multiple orders of magnitude but potentially problematic when denominators approach zero.
    \item Maximum hybrid error~\cite{xie_hyb_2024} $\epsilon_{\text{max-hyb}}(y, y') = \max_i \frac{|y_i - y'_i|}{1 + |y'_i|}$ measures the maximum difference between outputs by combining advantages of both absolute and relative errors and gracefully handling near-zero values. The maximum hybrid error delineates the maximum tolerance of the \texttt{allclose} function that is widely used in numerical libraries \cite{harris_array_2020,paszke_pytorch_2019}.
\end{enumerate}

\section{Statistical Analysis}

\subsection{Descriptive Statistics}

Before conducting formal hypothesis tests, we must thoroughly characterize each error distribution $\Delta_i$ through descriptive statistics that reveal its central tendency, spread, and shape. The mean $\bar{\Delta}_i = \frac{1}{N}\sum_{j=1}^N \Delta_i[j]$ provides our primary measure of average error magnitude, while the median $\tilde{\Delta}_i = \text{median}(\Delta_i)$ offers robustness against outliers and reveals whether the distribution exhibits skewness. Standard deviation $\sigma_i = \sqrt{\frac{1}{N-1}\sum_{j=1}^N (\Delta_i[j] - \bar{\Delta}_i)^2}$ quantifies variability, with high standard deviation indicating inconsistent behavior across different inputs---a potential red flag for reliability. Percentile analysis, particularly examining $P_k(\Delta_i)$ for $k \in \{50, 90, 95, 99\}$, illuminates the tail behavior that often matters most in practice: while average error may be acceptable, the 99th percentile error determines worst-case performance. Finally, the maximum error $\max(\Delta_i)$ captures the absolute worst-case scenario observed during testing.

By comparing these descriptive statistics across $\Delta_1$ and $\Delta_2$, we gain initial insights into whether one implementation consistently outperforms the other or whether their behaviors diverge in specific aspects (e.g., one may have lower average error but higher variability).

\subsection{Visualization}
Statistical summaries, while essential, cannot fully capture the richness of distributional information. Visualization techniques provide complementary insights that often reveal patterns invisible in summary statistics. Overlapping histograms display the complete distributional shapes, making it immediately apparent whether errors concentrate in narrow ranges or spread across wide spans, whether distributions exhibit skewness or multiple modes, and how much overlap exists between the two implementations. Box plots efficiently visualize quartiles and outliers, highlighting differences in central tendency, spread, and tail behavior in a compact format. Quantile-quantile (Q-Q) plots test whether the two distributions share similar shapes: if points fall along the diagonal, the distributions differ only in location and scale; deviations from the diagonal reveal more fundamental distributional differences. Scatter plots comparing $\Delta_1[i]$ versus $\Delta_2[i]$ for each input $i$ expose correlations in error behavior, revealing whether both implementations struggle with the same challenging inputs or exhibit independent error patterns.

\subsection{Hypothesis Testing}

If descriptive statistics and visualizations are insufficient to draw definitive conclusions, we can employ formal hypothesis tests to assess whether observed differences between $\Delta_1$ and $\Delta_2$ are statistically significant. For example, to test distributional equivalence (i.e., whether $\Delta_1$ and $\Delta_2$ arise from the same distribution), the two-sample Kolmogorov-Smirnov test provides a rigorous perspective. If the difference is significant, then we may wish to test the null hypothesis $H_0: \mu_1 \leq \mu_2$ against the alternative $H_a: \mu_1 > \mu_2$, where $\mu_i = \mathbb{E}[\Delta_i]$ represents the expected error of implementation $f_i$.  Since the error distributions in numerical computing are typically non-normal, we recommend the one-sided Wilcoxon signed-rank test or sign test as the default choice. This non-parametric test operates on the paired differences $D_i = \Delta_1[i] - \Delta_2[i]$ and tests whether these differences tend to be positive, without assuming any particular distributional form. When normality can be justified (e.g., via the Shapiro-Wilk test), a one-sided paired t-test offers greater statistical power.

\section{Practical Examples: Matrix Multiplication}

\subsection{Problem Setup}

Modern deep learning frameworks heavily rely on half-precision (FP16) matrix operations to accelerate neural network training and inference, yet these operations accumulate numerical errors that can impact model convergence and accuracy~\citep{micikevicius_mixed_2018}. We consider two real-world scenarios encountered when validating GPU implementations against CPU baselines.

In our experimental setup, we designate the GPU-based FP16 matrix multiplication as implementation $f_1$ (the custom implementation under test), the CPU-based FP16 matrix multiplication as implementation $f_2$ (the reference baseline), and FP64 matrix multiplication as our oracle $f_\Omega$. This configuration reflects a common validation scenario: developers implement optimized GPU kernels and must verify that numerical accuracy remains comparable to established CPU implementations while achieving substantial speedups.

\subsection{Case Study 1: Well-Behaved Configuration}

We first examine a well-behaved scenario using square $128 \times 128$ matrices, where both implementations should theoretically produce similar results.

\begin{lstlisting}[language=Python, caption=Dual-Delta Testing for Matrix Multiplication]
import torch
import numpy as np

def matmul_gpu_fp16(A, B):
    """Custom implementation: GPU FP16 matmul"""
    return torch.matmul(A, B)

def matmul_cpu_fp16(A, B):
    """Reference implementation: CPU FP16 matmul"""
    return torch.matmul(A.cpu(), B.cpu()).cuda()

def matmul_oracle(A, B):
    """Oracle: FP64 matmul"""
    return torch.matmul(A.double(), B.double())

def generate_input():
    """Generate random FP16 matrices"""
    A = torch.randn(128, 128).half().cuda()
    B = torch.randn(128, 128).half().cuda()
    return A, B

def get_error(res, res_oracle):
    """Compute maximum hybrid error"""
    diff = (res.double() - res_oracle).abs()
    normalized = diff / (1 + res_oracle.abs())
    return normalized.max().item()

def dual_delta_test(impl_1, impl_2, oracle, 
                   generate_input, get_error, num_tests):
    """Run Dual-Delta Testing"""
    delta_1 = []
    delta_2 = []
    
    for _ in range(num_tests):
        input_data = generate_input()
        res_1 = impl_1(*input_data)
        res_2 = impl_2(*input_data)
        res_oracle = oracle(*input_data)
        
        delta_1.append(get_error(res_1, res_oracle))
        delta_2.append(get_error(res_2, res_oracle))
    
    return delta_1, delta_2

# Run the test
delta_1, delta_2 = dual_delta_test(
    impl_1=matmul_gpu_fp16,
    impl_2=matmul_cpu_fp16,
    oracle=matmul_oracle,
    generate_input=generate_input,
    get_error=get_error,
    num_tests=1000
)

# Statistical analysis
mean_1, std_1 = np.mean(delta_1), np.std(delta_1)
mean_2, std_2 = np.mean(delta_2), np.std(delta_2)

print(f'GPU FP16: Mean Error = {mean_1:.6f}, Std = {std_1:.6f}')
print(f'CPU FP16: Mean Error = {mean_2:.6f}, Std = {std_2:.6f}')

# Visualization
_, bins = np.histogram(np.hstack((delta_1, delta_2)), bins=100)
plt.hist(delta_1, bins=bins.tolist(), alpha=0.5, label='GPU')
plt.hist(delta_2, bins=bins.tolist(), alpha=0.5, label='CPU')
plt.xlabel('Error'); plt.ylabel('Frequency')
plt.legend(); plt.title('Error Distribution Comparison')
plt.show()
\end{lstlisting}

\subsection{Results and Analysis}

Running this experiment with 1000 random matrix pairs yields results that validate both implementations. As shown in Figure~\ref{fig:good_case}, the error distributions of the GPU and CPU implementations exhibit nearly complete overlap, confirming that both achieve equivalent numerical accuracy for this configuration. The numerical results further corroborate this observation:

\begin{figure}[h]
    \centering
    \includegraphics[width=0.8\textwidth]{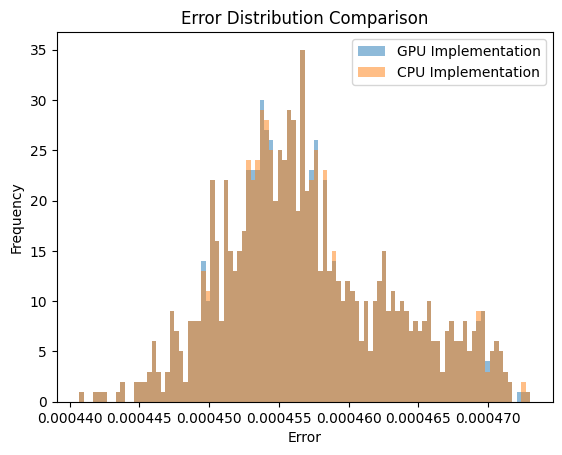}
    \caption{Error distribution comparison for $128 \times 128$ matrix multiplication. The GPU (blue) and CPU (orange) error distributions overlap almost entirely, indicating equivalent numerical accuracy.}
    \label{fig:good_case}
\end{figure}

\begin{table}[h]
    \centering
    \begin{tabular}{lcc}
        \toprule
        \textbf{Implementation} & \textbf{Mean Error} & \textbf{Std Dev} \\
        \midrule
        GPU (FP16) & $4.570838 \times 10^{-4}$ & $6.079639 \times 10^{-6}$ \\
        CPU (FP16) & $4.570844 \times 10^{-4}$ & $6.080522 \times 10^{-6}$ \\
        \bottomrule
    \end{tabular}
    \caption{Error statistics for $128 \times 128$ matrix multiplication.}
    \label{tab:good_case}
\end{table}

The mean errors are virtually identical (differing only in the sixth significant digit), and the standard deviations match to three significant figures. This near-perfect agreement confirms that neither implementation exhibits systematic accuracy advantages over the other for well-conditioned square matrix multiplication.

\subsection{Case Study 2: Detecting Numerical Issues}

A more revealing scenario emerges when we modify the matrix dimensions to test rectangular matrices of shape $128 \times 4096$ multiplied by $4096 \times 128$.

\begin{lstlisting}[language=Python, caption=Problematic Configuration]
def generate_input():
    A = torch.randn(128, 4096).half().cuda()
    B = torch.randn(4096, 128).half().cuda()
    return A, B

# Run dual-delta testing with new dimensions
delta_1, delta_2 = dual_delta_test(impl_1=matmul_gpu_fp16, 
    impl_2=matmul_cpu_fp16, oracle=matmul_oracle,
    generate_input=generate_input, get_error=get_error,
    num_tests=1000)
\end{lstlisting}

This seemingly innocuous change in matrix dimensions reveals a dramatic accuracy degradation in the GPU implementation. As shown in Table~\ref{tab:bad_case}, Dual-Delta Testing exposes that the default GPU implementation exhibits a mean error of $3.106403 \times 10^{-2}$, roughly $65 \times$ larger than the CPU reference's mean error of $4.768127 \times 10^{-4}$.

\begin{table}[h]
    \centering
    \begin{tabular}{lcc}
        \toprule
        \textbf{Implementation} & \textbf{Mean Error} & \textbf{Std Dev} \\
        \midrule
        GPU (FP16) & $3.106403 \times 10^{-2}$ & $5.460945 \times 10^{-3}$ \\
        CPU (FP16) & $4.768127 \times 10^{-4}$ & $3.262880 \times 10^{-6}$ \\
        \bottomrule
    \end{tabular}
    \caption{Error statistics for $128 \times 4096$ by $4096 \times 128$ matrix multiplication.}
    \label{tab:bad_case}
\end{table}

\begin{figure}[h]
    \centering
    \includegraphics[width=0.8\textwidth]{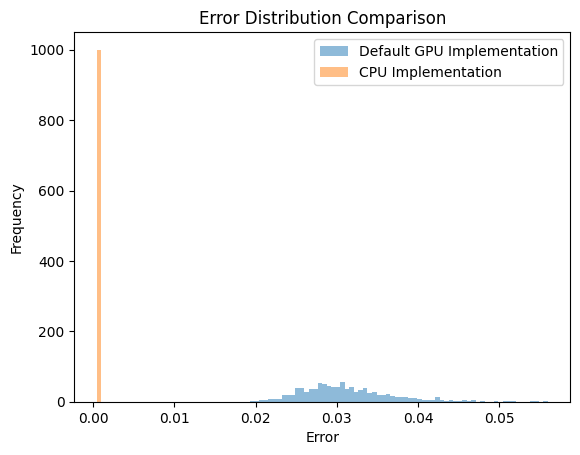}
    \caption{Error distribution comparison for $128 \times 4096$ by $4096 \times 128$ matrix multiplication. The GPU (blue) and CPU (orange) error distributions are clearly separated, revealing a significant accuracy degradation in the default GPU implementation.}
    \label{fig:bad_case}
\end{figure}

The error histograms (Figure~\ref{fig:bad_case}) show unambiguous separation between the two distributions, with the GPU distribution shifted far toward higher errors and exhibiting a much wider spread. This discrepancy points to accuracy-degrading optimizations in the GPU implementation---specifically, PyTorch's default \texttt{allow\_fp16\_reduced\_precision\_reduction} flag, which enables faster accumulation strategies at the cost of numerical precision. In the longer inner dimension ($K = 4096$), these reduced-precision reductions accumulate significantly more rounding error than in the well-behaved $128 \times 128$ case.

\subsection{Resolution and Validation}

Armed with these dual-delta results, we can diagnose and resolve the issue:

\begin{lstlisting}[language=Python, caption=Corrected Configuration]
# Disable reduced precision reduction
torch.backends.cuda.matmul.allow_fp16_reduced_precision_reduction = False

delta_1, delta_2 = dual_delta_test(impl_1=matmul_gpu_fp16, 
    impl_2=matmul_cpu_fp16, oracle=matmul_oracle,
    generate_input=generate_input, get_error=get_error,
    num_tests=1000)
\end{lstlisting}

\begin{figure}[h]
    \centering
    \includegraphics[width=0.8\textwidth]{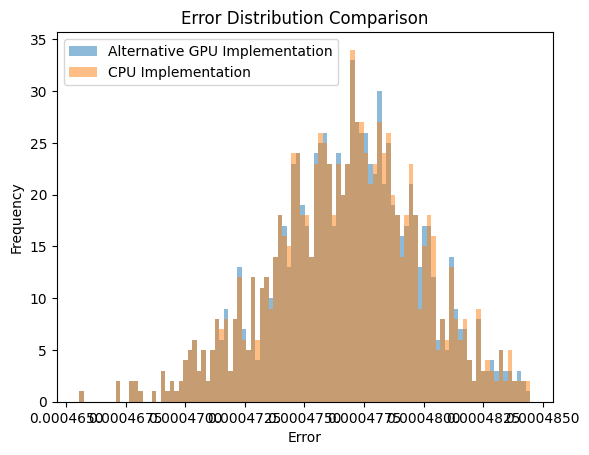}
    \caption{Error distribution comparison for $128 \times 4096$ by $4096 \times 128$ matrix multiplication after disabling reduced-precision reduction. The GPU (blue) and CPU (orange) error distributions now overlap almost entirely, confirming restored accuracy parity.}
    \label{fig:bad_case_fixed}
\end{figure}

After disabling the reduced-precision optimization, Dual-Delta Testing confirms that accuracy parity is restored. As shown in Figure~\ref{fig:bad_case_fixed} and Table~\ref{tab:bad_case_fixed}, both implementations now exhibit nearly identical error distributions, with mean errors returning to comparable magnitudes and standard deviations matching closely.

\begin{table}[h]
    \centering
    \begin{tabular}{lcc}
        \toprule
        \textbf{Implementation} & \textbf{Mean Error} & \textbf{Std Dev} \\
        \midrule
        GPU (FP16) & $4.767338 \times 10^{-4}$ & $3.124869 \times 10^{-6}$ \\
        CPU (FP16) & $4.767371 \times 10^{-4}$ & $3.127470 \times 10^{-6}$ \\
        \bottomrule
    \end{tabular}
    \caption{Error statistics for $128 \times 4096$ by $4096 \times 128$ matrix multiplication after disabling reduced-precision reduction.}
    \label{tab:bad_case_fixed}
\end{table}

This case study demonstrates the methodology's power: the dual-delta framework not only detected a subtle numerical issue that single-delta testing might miss, but also validated the proposed fix by confirming distributional equivalence.

\section{Conclusion}

Dual-Delta Testing provides a rigorous methodology for evaluating numerical accuracy in mixed-precision computing by measuring two error distributions against a high-precision oracle, overcoming the fundamental limitations of single-delta comparisons. Its dual-distribution design enables contextual evaluation that distinguishes precision-inherent errors from implementation-specific bugs, supports statistically sound hypothesis testing, and reveals detailed comparative insights into numerical accuracy. Our matrix multiplication case studies demonstrate its practical power: validating implementation equivalence in well-behaved configurations, detecting subtle accuracy regressions in challenging ones, and confirming proposed fixes---a diagnostic cycle that single-delta testing cannot reliably achieve. As mixed-precision computing expands across domains and hardware platforms, and as novel low-precision formats proliferate, we believe systematic validation methodologies like Dual-Delta Testing will evolve from occasional practice to standard expectation, much as unit testing has become standard in software engineering.

\bibliographystyle{plain}
\bibliography{references}

\end{document}